\documentclass[11pt]{amsart}
\usepackage{amssymb}
\usepackage{graphicx}
\usepackage{hyperref}
\usepackage{color}

\newtheorem{thm}{Theorem}[section]

\newtheorem{lem}[thm]{Lemma}

\theoremstyle{definition}
\newtheorem{defn}[thm]{Definition}
\theoremstyle{remark}
\newtheorem{rem}[thm]{Remark}

\numberwithin{equation}{section}

\def\br{{\mathbb R}}

\def\bz{{\mathbb Z}}

\def\p{\psi}

\def\bal1{\textbf{B}_1^+}

\allowdisplaybreaks

\usepackage{todonotes}

\begin{document}

\title[Anisotropic Herz-Amalgam spaces] {Atomic decomposition of anisotropic Herz-Amalgam spaces and boundedness of sublinear operators }

\author[M. Sultan, A. Hussain, B. Sultan and T. Mahmood]{Mehvish Sultan$^{1}$,   Amjad Hussain$^{1^*}$, Babar Sultan$^{1^*}$ and Tariq Mahmood$^{2}$  }%
\address{ $^{1^*}$Department of Mathematics, Quaid-I-Azam University, Islamabad  45320, Pakistan.}

\email{mehvishsultanbaz@gmail.com, a.hussain@qau.edu.pk,   babarsultan40@yahoo.com }%
\address{ $^{2}$
Department of Mathematics,
University of Chakwal.}
\email{tariq.mahmood@uoc.edu.pk}%



\date{\today}
\dedicatory{}%
------------------------------------------------------
\begin{abstract}
In this work, we introduce the idea of  anisotropic Herz-Amalgam  spaces. Then we find the atomic decomposition in these spaces.  As an application, in these spaces, we demonstrate the boundedness of certain sublinear operators. 
\end{abstract}
\maketitle
\section{\textbf{Introduction}}

In 1926, Wiener \cite{wiener} introduced the idea of amalgam spaces. Let $p,q\in (0,\infty)$, then Amalgam space $\left(L^p,L^q \right)\left(\mathbb{R} \right)$ is defined as 
$$\left(L^p,L^q \right)\left(\mathbb{R} \right):=\left\{ f \in L^p_{\mathrm{loc}}: \left[ \sum_{m\in \mathbb{Z}}\left( \int_m^{m+1}\left|f(x)\right|^{p}dx\right)^{q/p}\right]^{1/q}<\infty\right\}. $$

Many writers have studied amalgam spaces or some of its applications \cite{bertrandias,cowling,fournier}. The fact that amalgam spaces provide information regarding the local $L^p$ and global $L^q$ features of the functions, whereas $L^p$ spaces do not, is a key distinction when comparing them to $L^p$ spaces.

A locally integrable function on $\mathbb{R}^{n}$, is a  weight $\Omega $ which  takes values almost everywhere in $(0, \infty)$. Let $1 \leqslant p<\infty$, then weighted Lebesgue spaces 
$L_\Omega  ^{p}\left(\mathbb{R}^{n}\right)$ the space of all functions such that the norm 
$
\left\| g\right\|_{L_\Omega  ^{p}\left(\mathbb{R}^{n}\right)}=\left(\int_{\mathbb{R}^{n}}|g(v)|^{p} \Omega (y) d y\right)^{1 / p}$ is finite.

Let $1<p<\infty$, $1\leq q \leq \infty$, $0<t<\infty$ and $v,w$ are two weights. Then the weighted Amalgam space $\left(L^p_w,L^q_v \right)_t\left(\mathbb{R}^n \right)$ is the space of all measurable functions $g$ such that 

$$\left\|f\right\|_{\left(L^p_w,L^q_v \right)_t\left(\mathbb{R}^n \right)}:=\left\|\left( \frac{1}{w(B(\cdot,t))}\int\limits_{B(\cdot,t)}\left|f(y)\right|^pw(y)dy\right)^{1/p} \right\|_{L^q_v(\mathbb{R}^n)}<\infty,$$

where  $B(y, r):=\left\{x \in \mathbb{R}^{n}:|y-x|<r\right\}$.

The recent advancements in variable exponent analysis were initiated in \cite{23}, who defined a class of function spaces with variable exponents and established the fundamental properties of Lebesgue spaces with variable exponents. In particular, the theory of variable exponent analysis has focused on the boundedness of the Hardy–Littlewood maximal operator 
$M$ in these spaces \cite{4,5,7}. Muckenhoupt initially established the classical $A_{p}$ weight theory in \cite{32} while studying weighted $L^{p}$ boundedness of Hardy-Littlewood maximal functions. Herz \cite{herz} created a class of function spaces known as the Herz spaces in 1968 while studying completely convergent Fourier transformations. For more results  generalized versions of Herz space like boundedness of sublinear operators, fractional integrals, the commutator of singular integrals with BMO functions, and the commutator of fractional integrals with BMO functions see \cite{aims4,ghm4,ghm5,ghm2,babargh2,babargh4,ukr,krag,ineq}.

Let  $\alpha \in \br $, $q \in (0, \infty]$, and $0 < p\leq \infty $,   then the homogeneous version of   Herz spaces $\dot{K}^{\alpha,q }_{ p}$ are defined by
	$$\dot{K}^{\alpha,q }_{ p}=\left\{ g \in L^{p} _{\mathrm{loc}}(\br^n\setminus\{0\}): \| g \| _ {\dot{K}^{\alpha,q }_{ p}} < \infty \right\},$$ 
	where 
	$$\| g \| _ {\dot{K}^{\alpha,q }_{ p}}= \left(   \sum \limits _{\ell= -\infty}^{\infty} 2^{\ell \alpha q} \| g\chi _\ell \| ^{q}_{L^p}\right)^{\frac{1}{q}}.$$

In this article, we introduce the idea of anisotropic Herz-Amalgam spaces  by using anisotropic Herz spaces and Amalgam spaces. We will establish the atomic decomposition in these spaces. Then, we obtain  boundedness for sublinear  in these spaces.

\section{\textbf{Preliminaries}}

Let $\ell \in \bz$, $R_\ell: B_\ell\setminus B_{\ell-1}$ where $B_\ell: \left\{ x \in \br^n \; :\; |x| \leq 2^\ell\right\}$. $\chi_\ell := \chi_{R_\ell}$, where $\chi_{R_\ell}$ is the characteristic function of $R_\ell$. $C$ is the constant which vary from the line to line.

A $n\times n$ real matrix $A$ is called an expansive matrix or dilation, if all eignevalues $\lambda$ of $A$ satisfy $|\lambda|>1$. Let $\lambda_1,\cdots,\lambda_n$ are eigenvalues of $A$ such that $1<|\lambda_1|\leq \cdots \leq |\lambda_n|$. A set $\Delta \subset \mathbb{R}^n$ is called ellipsoid if $\Delta = \left\{x]\in \mathbb{R}^n : |Px|<1\right\}$ where $P$ is nondegenerate $n\times n$ matrix and $|\cdot|$ is Euclidean norm in $\mathbb{R}^n$. For dilation $A$, there is an ellipsoid $\Delta$ and $r>1$ so that $\Delta \subset r\Delta \subset A\Delta$ and $|\Delta|=1$, $|\Delta|$ is the Lebesgue measure of $\Delta$. 
If $B_k=A^k \Delta$ for $k\in \mathbb{Z}$, then $b_k \subset rB_k \subset B_{k+1}$, and $|B_k|=b^k$ such that $b=|\mathrm{det}A|>1$. Assume that $w$ is the smallest integer such that $2B_0 \subset A^w B_0=B_w$. A quaisi-norm with expansive matrix $A$ is a measurable mapping $\rho_A:\mathbb{R}^n \rightarrow [0,\infty)$ such that
$$\rho_A(y)>0\;\;\; \it{for} \; x\neq 0,$$
$$\rho_A(Ay)= |\mathrm{det}A|\rho(y) \;\;\; \it{for} \; x\in \mathbb{R}^n,$$
$$\rho_A(x+y)\leq C\left(\rho_A(x)+\rho_A(y)\right)\;\;\; \it{for} \; x,y\in \mathbb{R}^n,$$
where $C\geq 1.$ To define step homogeneous quasi-norm $\rho$ on $\mathbb{R}^n$ induced by dilation $A$ as
$$\left(\rho(x)=b^j\;\; \it{for} \; x\in B_{j+1}\setminus B_j\right)\;\;\; \it{and} \;\;\; \left(0 \;\;\; \it{for}\; x=0\right).$$
If $x,y \in \mathbb{R}^n,$ then we have
\begin{equation}\label{1.1}
    \rho(x+y)\leq b^w \left(\rho(x)+\rho(y)\right).
\end{equation}

\begin{defn}
Let  $\alpha \in \br $,  $1<p<\infty$, $1\leq q \leq \infty$, $0<t,r<\infty$ and $v,w$ are two weights.   Then the homogeneous version of  weighted Herz-Amalgam spaces $\dot{K}{\left(L^p_w,L^q_v \right)_t\left(\mathbb{R}^n \right)}$ are defined by
	$$\dot{K}{\left(L^p_w,L^q_v \right)_t\left(\mathbb{R}^n \right)}=\left\{ f \in \left(L^p_w,L^q_v \right)_t\left(\mathbb{R}^n \right): \| f \| _ {\dot{K}{\left(L^p_w,L^q_v \right)_t\left(\mathbb{R}^n \right)}} < \infty \right\},$$ 
	where 
	$$\| f \| _ {\dot{K}{\left(L^p_w,L^q_v \right)_t\left(\mathbb{R}^n \right)}}= \left(   \sum \limits _{\ell= -\infty}^{\infty} b^{\ell \alpha r} \| f \chi _\ell \| ^{r}_{\left(L^p_w,L^q_v \right)_t\left(\mathbb{R}^n \right)}\right)^{\frac{1}{r}}.$$

\end{defn}

Let $1\leq p,q <\infty$. The Hölder inequality for  weighted Amalgam space \cite{archive} is given as:

     $$
\int_{\mathbb{R}^{n}}|f(z) g(z)| dz \leq C\left\| f\right\| _{\left(L^p_w,L^q_v \right)_t\left(\mathbb{R}^n \right)}\|g\|_{\left(L^{p^\prime}_{w^\prime},L^{q^\prime}_{v^\prime} \right)_t\left(\mathbb{R}^n \right)}
$$

where $\frac{1}{p}=\frac{1}{p^\prime}=\frac{1}{q}=\frac{1}{q^\prime}=1,v^\prime=v^{1-p^\prime}, w^\prime=v^{1-q^\prime}$.

\begin{lem}\label{lem1}
Let $1<p<\infty$, $1\leq q \leq \infty$, $0<t<\infty$ and $v,w$ are two weights. A characteristic function on $B\left(y_{0}, s_{0}\right)$ fulfills

  $$
\left\|\chi_{B\left(y_{0}, s_{0}\right)}\right\|_{\left(L^p_w,L^q_v \right)_t\left(\mathbb{R}^n \right)} \leq C s_{0}^{n / q}.
$$

\begin{proof}
   Let $x \in \mathbb{R}^{n}$, we get

$$
\begin{aligned}
\left\|\chi _{B\left(y_{0}, s_{0}\right)}\right\|_{\left(L^p_w,L^q_v \right)_t\left(\mathbb{R}^n \right)} & \left.=\left\{\int_{\mathbb{R}^{n}} \frac{\left\|\chi _{B\left(y_{0}, s_{0}\right)} \chi _{B(x,t)}\right\|_{L^p_w }}{\left\|\chi _{B(x,t)}\right\|_{L^p_w }}\right]^{q}v(x) d x\right\}^{\frac{1}{q}} \\
& =\frac{1}{\left\|\chi _{B(\overrightarrow{0}, t)}\right\|_{L^p_w }}\left[\int_{\mathbb{R}^{n}}\left\|\chi _{B\left(\overrightarrow{0}, s_{0}\right)} \chi _{B\left(x-y_{0}, t\right)}\right\|_{L^p_w }^{q}v(x) d x\right]^{\frac{1}{q}} \\
& =\frac{1}{\left\|\chi _{B(\overrightarrow{0}, t)}\right\|_{L^p_w }}\left[\int_{\mathbb{R}^{n}}\left\|\chi _{B\left(\overrightarrow{0}, s_{0}\right)} \chi _{B(x,t)}\right\|_{L^p_w }^{q} v(x)d x\right]^{\frac{1}{q}} \\
& =\left\|\chi _{B\left(\overrightarrow{0}, s_{0}\right)}\right\|_{\left(L^p_w,L^q_v \right)_t\left(\mathbb{R}^n \right)}
\end{aligned}
$$

This means that, assuming nothing is lost in general, $B_{1}:=B(\overrightarrow{0}, 1)$ and $B_{2}:=B\left(s_{0},\overrightarrow{0}\right)$ with $ 1<s_{0}<\infty$ are used. The geometric property tells us that there exists $M \in \mathbb{N}$ such that $B\left(\overrightarrow{0}, s_{0}\right) \subset \bigcup\limits_{i=1}^{M} B\left(x_{i}, 1\right)$, and that $M \sim\left|B_{2}\right|^{n}$ and $\left\{x_{1}, \ldots, x_{M}\right\}$.

$$
\begin{aligned}
\left\|\chi _{B\left(\overrightarrow{0}, s_{0}\right)}\right\|_{\left(L^p_w,L^q_v \right)_t\left(\mathbb{R}^n \right)} & =\left\|\chi _{B_{2}}\right\|_{\left(L^{p/q}_w,L^1_v \right)_t\left(\mathbb{R}^n \right)}^{\frac{1}{q}} \\
& \leq\left\|\sum_{i=1}^{M} \chi _{B\left(x_{i}, 1\right)}\right\|_{\left(L^{p/q}_w,L^1_v \right)_t\left(\mathbb{R}^n \right)}^{\frac{1}{q}} \lesssim\left(\sum_{i=1}^{M}\left\|\chi _{B\left(x_{i}, 1\right)}\right\|_{\left(L^{p/q}_w,L^1_v \right)_t\left(\mathbb{R}^n \right)}\right)^{\frac{1}{q}} \\
& \sim\left|B_{2}\right|^{\frac{1}{q}}\left\|\chi _{B(\overrightarrow{0}, 1)}\right\|_{\left(L^p_w,L^q_v \right)_t\left(\mathbb{R}^n \right)} \sim s_{0}^{n / q} .
\end{aligned}
$$

This completes the proof.
\end{proof}
\end{lem}

\begin{rem}
    Let $1<p<\infty$, $1\leq q \leq \infty$, $0<t<\infty$, $v,w$ are two weights and $ \ell \in  \mathbb{Z}$, a characteristic function on $R_{\ell}$ satisfies

$$
\left\|\chi_{R_{\ell}}\right\|_{\left(L^p_w,L^q_v \right)_t\left(\mathbb{R}^n \right)} \leq\left\|\chi_{B_{\ell}}\right\|_{\left(L^p_w,L^q_v \right)_t\left(\mathbb{R}^n \right)} \leq C b^{\ell n  / q}.
$$
\end{rem}

\begin{lem}
 Let $1<p<\infty$, $1\leq q \leq \infty$, $0<t<\infty$, $v,w$ are two weights.   For any ball $B$ in $\mathbb{R}^{n}$, there exists a positive constant $C$ such that for $\left\|\left| f\right|_B\chi_B \right\|_{\left(L^p_w,L^q_v \right)_t}\leq \left\|(f\chi_B) \right\|_{\left(L^p_w,L^q_v \right)_t}$ the inequality

$$\frac{1}{ \vert B_k \vert } \Vert \chi_{B_k} \Vert _{\left(L^p_w,L^q_v \right)_t} \Vert \chi_{B_k} \Vert _{\left(L^{p^\prime}_{w^\prime},L^{q^\prime}_{v^\prime} \right)_t}\leq C$$
holds.

\begin{proof}
Let $B=B_{k_0}(x_0)$  $\forall \; x \in B$, we get $B=B_{k_0}(x)$  $\forall \; x \in B$. Now for all $B_k(a)$, we get 

\begin{align*}
    \frac{1}{|B_k(a)|} &\left\| \chi_{B_k}(a)\right\||_{\left(L^p_w,L^q_v \right)_t}\left\| \chi_{B_k}(a)\right\||_{\left(L^{p^\prime}_{w^\prime},L^{q^\prime}_{v^\prime} \right)_t}\\
    \leq & \frac{1}{|B_k(a)|} \left\| \chi_{B_k}(a)\right\||_{\left(L^p_w,L^q_v \right)_t} \sup \left\{ \int_{\mathbb{R}^n}\left| f(x) \chi_B(a)\right| dx\; : \left\|f \right\|_{\left(L^p_w,L^q_v \right)_t}\leq 1\right\}\\
    = & \left\| \chi_{B_k}(a)\right\||_{\left(L^p_w,L^q_v \right)_t} \sup \left\{ \left\|\left| f\right|{\chi_B(a)}\right\|_{\left(L^p_w,L^q_v \right)_t} \; : \left\|f \right\|_{\left(L^p_w,L^q_v \right)_t}\leq 1\right\}\\
    =& \sup \left\{ \left\|\left| f\right|{\chi_B(a)}\right\|_{\left(L^p_w,L^q_v \right)_t} \; : \left\|f \right\|_{\left(L^p_w,L^q_v \right)_t}\leq 1\right\}\\
     \leq & \sup \left\{ \left\| f {\chi_B(a)}\right\|_{\left(L^p_w,L^q_v \right)_t} \; : \left\|f \right\|_{\left(L^p_w,L^q_v \right)_t}\leq 1\right\}\\
     \leq 1 .
\end{align*}
    Hence we obtain our desired result.
\end{proof}
\end{lem}

\section{Atomic Decomposition of anisotropic Herz-Amalgam   spaces}
\begin{defn}
 Let $\alpha\in \mathbb{R}$, $1<p<\infty$, $1\leq q \leq \infty$, $0<t<\infty$, $v,w$ are two weights.   
 \begin{enumerate}
       \item [(i)] A measurable function $a(y)$ is called central $\left(\alpha,{\left(L^p_w,L^q_v \right)_t}\right)$-block if $\mathrm{supp} a\subset B_l$ and $\left\|a\right\|_{\left(L^p_w,L^q_v \right)_t}\leq b^{-l \alpha}.$
       \item[(ii)] A measurable function $a(y)$ is called central $\left(\alpha,{\left(L^p_w,L^q_v \right)_t}\right)$-block if $\mathrm{supp} a\subset B_l$.
   \end{enumerate}
\end{defn}

\begin{thm}

 Let $\alpha\in \mathbb{R}$, $1<p<\infty$, $1\leq q \leq \infty$, $0<t,r<\infty$, $v,w$ are two weights.  The following two statements are equivalent:\\
(i) $g \in {\dot{K}{\left(L^p_w,L^q_v \right)_t\left(\mathbb{R}^n \right)}}$.\\
(ii) $g$ are given as

\begin{equation*}
g(y)=\sum_{k=-\infty}^{\infty} \lambda_{k} b_{k}(y) \tag{2.1}
\end{equation*}

where each $b_{k}$ is a central $\left(\alpha,{\left(L^p_w,L^q_v \right)_t}\right)$-block with support contained in $B_{k}$ and $\sum_{k=-\infty}^{\infty}\left|\lambda_{k}\right|^{r}<\infty$.
\begin{proof}
    
 We first prove $(i)$ implies $(ii)$. For every $g \in {\dot{K}{\left(L^p_w,L^q_v \right)_t\left(\mathbb{R}^n \right)}}$, write

$$
\begin{aligned}
g(y) & =\sum_{k=-\infty}^{\infty} g(y) \chi_{k}(y) \\
& =\sum_{k=-\infty}^{\infty}b^{k \alpha }\left\| g\chi _{k}\right\|_{\left(L^p_w,L^q_v \right)_t} \frac{g(y) \chi_{k}(y)}{b^{k \alpha }\left\| \chi_{k}\right\|_{{\left(L^p_w,L^q_v \right)_t}}} \\
& =\sum_{k=-\infty}^{\infty} \lambda_{k} b_{k}(x),
\end{aligned}
$$

where $\lambda_{k}=b^{k \alpha }\left\| g\chi _{k}\right\|_{\left(L^p_w,L^q_v \right)_t}$ and $b_{k}(x)=\frac{g(y) \chi_{k}(y)}{b^{k \alpha }\left\| \chi _{k}\right\|_{\left(L^p_w,L^q_v \right)_t}}$.\\
It is obvious that supp $b_{k} \subset B_{k}$ and $\left\|b_{k}\right\|_{\left(L^p_w,L^q_v \right)_t}=\left|B_{k}\right|^{-\alpha / n}$. Thus, each $b_{k}$ is a central $\left(\alpha,{\left(L^p_w,L^q_v \right)_t}\right)$-block with the support $B_{k}$ and

$$
\sum_{k=-\infty}^{\infty}\left|\lambda_{k}\right|^{r}=\sum_{k=-\infty}^{\infty}b^{k \alpha }\left\| g\chi _{k}\right\|_{\left(L^p_w,L^q_v \right)_t}^{r}=\left\|g\right\|^{r}_{\dot{K}{\left(L^p_w,L^q_v \right)_t\left(\mathbb{R}^n \right)}}<\infty.
$$

Now we prove $(ii)$ implies $(i)$. Let $g(y)=\sum_{k=-\infty}^{\infty} \lambda_{k} b_{k}(y)$, we get

\begin{equation}\label{2.2}
\left\|g\chi _{j}\right\|_{\left(L^p_w,L^q_v \right)_t} \leq \sum_{k=j}^{\infty}\left|\lambda_{k}\right|\left\|b_{k}\right\|_{\left(L^p_w,L^q_v \right)_t}. 
\end{equation}

If $0<r \leq 1$. From \eqref{2.2} it follows that
\begin{align*}
    \left\|g\right\|_{\dot{K}{\left(L^p_w,L^q_v \right)_t\left(\mathbb{R}^n \right)}}^{r} & =\sum\limits_{k=-\infty}^{\infty}b^{k \alpha }\left\| g\chi _{k}\right\|_{\left(L^p_w,L^q_v \right)_t}^{r} \\
    & \leq \sum\limits_{k=-\infty}^{\infty}b^{k \alpha r}\left\| g\chi _{k}\right\|_{\left(L^p_w,L^q_v \right)_t}^{r} \\
        & \leq \sum\limits_{k=-\infty}^{\infty}b^{k \alpha r}\left(\sum_{j=k}^{\infty}\left|\lambda_{j}\right|^{r}\left\|b_{j}\right\|^{r}_{\left(L^p_w,L^q_v \right)_t} \right).
\end{align*}

Let $I= \sum\limits_{k=-\infty}^{\infty}b^{k \alpha r}\left(\sum_{j=k}^{\infty}\left|\lambda_{j}\right|^{r}\left\|b_{j}\right\|^{r}_{\left(L^p_w,L^q_v \right)_t} \right).$ By using the fact  $0<\alpha<\infty$, we get

\begin{align*}
I & =\sum\limits_{k=-\infty}^{\infty}b^{k \alpha r}\left(\sum_{j=k}^{\infty}\left|\lambda_{j}\right|^{r}\left\|b_{j}\right\|^{r}_{\left(L^p_w,L^q_v \right)_t} \right)\\
& \lesssim  \sum_{k=-\infty}^{\infty} b^{k \alpha r} \sum_{j=k}^{\infty}\left|\lambda_{j}\right|^{r} b^{-j \alpha r} \\
& \lesssim  \sum_{k=-\infty}^{\infty} \sum_{j=k}^{\infty}\left|\lambda_{j}\right|^{r} b^{(k-j) \alpha r} \\
& \lesssim \sum_{j=-\infty}^{\infty} \sum_{k=-\infty}^{j}\left|\lambda_{j}\right|^{r} b^{(k-j) \alpha r}\\
& \lesssim \sum_{j=-\infty}^{\infty}\left|\lambda_{j}\right|^{r} .
\end{align*}

If $1<r<\infty$, we have

For $I$, by \eqref{2.2}, $0<\alpha <\infty$ and the Hölder inequality, we get

\begin{align*}
 I &\lesssim \sum_{k=-\infty}^{\infty} b^{k \alpha r}\left(\sum_{j=k}^{\infty}\left|\lambda_{j}\right|\left\|b_{j}\right\|_{\left(L^p_w,L^q_v \right)_t}\right) \\
& \lesssim \sum_{k=-\infty}^{\infty}\left(\sum_{j=k}^{\infty}\left|\lambda_{j}\right| b^{(k-j) \alpha }\right)^{r} \\
& \lesssim  \sum_{k=-\infty}^{\infty}\left(\sum_{j=k}^{\infty}\left|\lambda_{j}\right|^{r} b^{(k-j) \alpha r / 2}\right)\left(\sum_{j=k}^{\infty} b^{(k-j) \alpha r^{\prime} / 2}\right)^{r/r^{\prime}} \\
& \lesssim  \sum_{k=-\infty}^{\infty}\left(\sum_{j=k}^{\infty}\left|\lambda_{j}\right|^{r} b^{(k-j) \alpha r / 2}\right) \\
& \lesssim \sum_{j=-\infty}^{\infty} \sum_{k=-\infty}^{j}\left|\lambda_{j}\right|^{r} b^{(k-j) \alpha r / 2} \\
& \lesssim \sum_{j=-\infty}^{\infty}\left|\lambda_{j}\right|^{r} .
\end{align*}

Thus the proof of the Theorem is completed. Non-homogeneous version of the Theorem is obtained similarly.  
\end{proof}

\end{thm}

\begin{rem}
    From the proof of the above Theorem, it is easy to see that if $g \in {\dot{K}{\left(L^p_w,L^q_v \right)_t\left(\mathbb{R}^n \right)}}$ and $g(y)=\sum_{k=-\infty}^{\infty} \lambda_{k} b_{k}(x)$ be a central $\left(\alpha,{\left(L^p_w,L^q_v \right)_t}\right)$-block decomposition, then

$$
\left\|g\right\|_{\dot{K}{\left(L^p_w,L^q_v \right)_t\left(\mathbb{R}^n \right)}} \approx \left(\sum_{k=-\infty}^{\infty}\left|\lambda_{k}\right|^r\right)^{1/r}.
$$
\end{rem}

\section{ Boundedness of sublinear operators on anisotropic Herz-Amalgam   spaces}
As applications of the decomposition theorems, we will prove the boundedness of sublinear operators on anisotropic Herz-Amalgam   spaces. 
\begin{thm}
  Let $\alpha\in \mathbb{R}, $ with $0<\alpha<1/q^\prime$ and $ p,q,r \in (0, \infty]$.   If a sublinear operator $T$ satisfies
    \begin{equation}\label{3.1}
|T g(y)| \lesssim \int_{\mathbb{R}^{n}} \frac{|g(y)|}{\rho(x-y)} d y, \quad x \notin \operatorname{supp} g ,
\end{equation}

for any $g \in \left(L^p_w,L^q_v \right)_t $ with a compact support and $T$ is bounded on $\left(L^p_w,L^q_v \right)_t $.  Then $T$ is bounded in ${\dot{K}{\left(L^p_w,L^q_v \right)_t\left(\mathbb{R}^n \right)}}$.

\begin{proof}
   Let $g \in {\dot{K}{\left(L^p_w,L^q_v \right)_t\left(\mathbb{R}^n \right)}}$. By using the atomic decomposition Theorem, $g(y)=$ $\sum_{j=-\infty}^{\infty} \lambda_{j} b_{j}(x)$, where each $b_{j}$ is a central $(\alpha, \p$-block with support contained in $B_{j}$ and

$$
\left\|g\right\|_{\dot{K}{\left(L^p_w,L^q_v \right)_t\left(\mathbb{R}^n \right)}} \approx\left(\sum_{j=-\infty}^{\infty}\left|\lambda_{j}\right|^{r}\right)^{1/r}
$$

Therefore, we get

\begin{align*}
& \|Tg\|_{\dot{K}{\left(L^p_w,L^q_v \right)_t\left(\mathbb{R}^n \right)}}^{r}=\sum_{k=-\infty}^{\infty}b^{k \alpha r}\left\|(Tg) \chi_{k}\right\|_{\left(L^p_w,L^q_v \right)_t}^{r} \\
& \lesssim  \sum_{k=-\infty}^{\infty} b^{k \alpha r}\left\|(Tg) \chi_{k}\right\|_{\left(L^p_w,L^q_v \right)_t }^{r} \\
& \lesssim \sum_{k=-\infty}^{\infty} b^{k \alpha r}\left(\sum_{j=-\infty}^{k-w-1}\left|\lambda_{j}\right|\left\|\left(T b_{j}\right) \chi_{k}\right\|_{\left(L^p_w,L^q_v \right)_t }\right)^{r} \\
& +\sum_{k=-\infty}^{\infty} b^{k \alpha r}\left(\sum_{j=k-w}^{\infty} \left| \lambda_{j}\right| \left\|\left(T b_{j}\right) \chi_{k} \right\|_{\left(L^p_w,L^q_v \right)_t }\right)^{r} \\
& =I_{1}+I_{2} .
\end{align*}

Let us first estimate $I_{1}$. If $j \leq k-w-1, x \in C_{k}$ and $y \in B_{j}$, by \eqref{1.1} we get

$$
\rho(x-y) \geq b^{-w} \rho(x)-\rho(y) \geq b^{-w} \rho(x)-b^{-w-1} \rho(x)=b^{-w}(1-1 / b) \rho(x).
$$

Therefore by \eqref{3.1} and the generalized Hölder inequality, we have

\begin{align*}
\left|T b_{j}(x)\right| & \leq C \rho(x)^{-1} \int_{B_{j}}\left|b_{j}(y)\right| d y \\
& \leq C b^{-k}\left\|b_{j}\right\|_{\left(L^p_w,L^q_v \right)_t }\left\|\chi_{B_{j}}\right\|_{\left(L^{p^\prime}_{w^\prime},L^{q^\prime}_{v^\prime} \right)_t}.
\end{align*}

Applying  Lemmas $2.3$ and $2.4$, we get

\begin{align*}
 \left\|\left(T b_{j}\right) \chi_{k}\right\|_{\left(L^p_w,L^q_v \right)_t } &\lesssim b^{-k}\left\|b_{j}\right\|_{\left(L^p_w,L^q_v \right)_t }\left\|\chi_{B_{j}}\right\|_{\left(L^{p^\prime}_{w^\prime},L^{q^\prime}_{v^\prime} \right)_t}\left\|\chi_{B_{k}}\right\|_{\left(L^p_w,L^q_v \right)_t } \\
& \lesssim b^{-k}\left\|b_{j}\right\|_{\left(L^p_w,L^q_v \right)_t }\left(\left|B_{k}\right|\left\|\chi_{B_{k}}\right\|_{\left(L^{p^\prime}_{w^\prime},L^{q^\prime}_{v^\prime} \right)_t}^{-1}\right)\left\|\chi_{B_{j}}\right\|_{\left(L^{p^\prime}_{w^\prime},L^{q^\prime}_{v^\prime} \right)_t} \\
& \lesssim\left\|b_{j}\right\|_{\left(L^p_w,L^q_v \right)_t } \frac{\left\|\chi_{B_{j}}\right\|_{\left(L^{p^\prime}_{w^\prime},L^{q^\prime}_{v^\prime} \right)_t}}{\left\|\chi_{B_{k}}\right\|_{\left(L^{p^\prime}_{w^\prime},L^{q^\prime}_{v^\prime} \right)_t}} \\
& \lesssim b^{(j-k)/q^\prime}\left\|b_{j}\right\|_{\left(L^p_w,L^q_v \right)_t } . 
\end{align*}

Therefore, for $0<r \leq 1$, and $0<\alpha <1/r^\prime$, we have

\begin{align*}
I_{1} & =\sum_{k=-\infty}^{\infty} b^{k \alpha r}\left(\sum_{j=-\infty}^{k-w-1}\left|\lambda_{j}\right|\left\|\left(T b_{j}\right) \chi_{k}\right\|_{\left(L^p_w,L^q_v \right)_t }\right)^{r} \\
& \lesssim\ \sum_{k=-\infty}^{\infty} b^{k \alpha r}\left(\sum_{j=-\infty}^{k-w-1}\left|\lambda_{j}\right|^{r} b^{\left[(j-k)/q^\prime-j \alpha \right]r}\right) \\
& \lesssim\ \sum_{j=-\infty}^{-w-2}\left|\lambda_{j}\right|^{r} \sum_{k=j+w+1}^{\infty} b^{(j-k)\left[1/q^\prime-\alpha \right]r} \\
& \lesssim\ \sum_{j=-\infty}^{-w-2}\left|\lambda_{j}\right|^{r}\\
&\lesssim\|g\|_{\dot{K}{\left(L^p_w,L^q_v \right)_t\left(\mathbb{R}^n \right)}}^{r}.
\end{align*}

When $1<r<\infty$,  $0<\alpha <1/q^\prime$, and using the Hölder inequality, we get

\begin{align*}
I_{1} & \lesssim \sum_{k=-\infty}^{\infty} b^{k \alpha r}\left(\sum_{j=-\infty}^{k-w-1}\left|\lambda_{j}\right| b^{(j-k)/q^\prime-j \alpha }\right)^{r} \\
& \lesssim \sum_{k=-\infty}^{\infty}\left(\sum_{j=-\infty}^{k-w-1}\left|\lambda_{j}\right|^{r} b^{(j-k)\left[1/q^\prime-\alpha \right]r/ 2}\right)\left(\sum_{j=-\infty}^{k-w-1} b^{(j-k)\left[1/q^\prime-\alpha \right]r^{\prime} / 2}\right)^{r/r^{\prime}} \\
& \lesssim \sum_{k=-\infty}^{\infty}\left(\sum_{j=-\infty}^{k-w-1}\left|\lambda_{j}\right|^{r} b^{(j-k)\left[1/q^\prime-\alpha \right]r / 2}\right) \\
& \lesssim \sum_{j=-\infty}^{-w-2}\left|\lambda_{j}\right|^{r} \sum_{k=j+w+1}^{-1} b^{(j-k)\left[1/q^\prime-\alpha \right]r/ 2} \\
& \lesssim \sum_{j=-\infty}^{-w-2}\left|\lambda_{j}\right|^{r} \\
&\lesssim\|g\|_{\dot{K}{\left(L^p_w,L^q_v \right)_t\left(\mathbb{R}^n \right)}}^{r}.
\end{align*}

Next we find the estimate of  $I_{2}$. If  $0<r \leq 1$, by $\left(L^p_w,L^q_v \right)_t$ boundedness of $T$, we get
\begin{align*}
I_{2} & =\sum_{k=-\infty}^{\infty} b^{k \alpha r}\left(\sum_{j=k-w}^{\infty}\left|\lambda_{j}\right|\left\|\left(T b_{j}\right) \chi_{k}\right\|_{\left(L^p_w,L^q_v \right)_t }\right)^{r} \\
& \lesssim \sum_{k=-\infty}^{\infty} b^{k \alpha r}\left(\sum_{j=k-w}^{\infty}\left|\lambda_{j}\right|^{r} \|\left. b_{j}\right|_{\left(L^p_w,L^q_v \right)_t } ^{r} \right) \\
& \lesssim \sum_{k=-\infty}^{\infty} b^{k \alpha r}\left(\sum_{j=k-w}^{\infty}\left|\lambda_{j}\right|^{r} b^{-j \alpha r}\right) \\
& \lesssim \sum_{k=-\infty}^{\infty} \sum_{j=k-w}^{\infty}\left|\lambda_{j}\right|^{r} b^{(k-j) \alpha r}\\
& \lesssim \sum_{j=-\infty}^{\infty}\left|\lambda_{j}\right|^{r} \sum_{k=-\infty}^{j+w} b^{(k-j) \alpha r}\\
& \lesssim \sum_{j=-\infty}^{\infty}\left|\lambda_{j}\right|^{r} \\
&\lesssim\|g\|_{\dot{K}{\left(L^p_w,L^q_v \right)_t\left(\mathbb{R}^n \right)}}^{r}.
\end{align*}

If $1<r<\infty$,  by using  $\left(L^p_w,L^q_v \right)_t $ boundedness of $T$ and the Hölder's inequality, we get

\begin{align*}
I_{2} & \lesssim \sum_{k=-\infty}^{\infty} b^{k \alpha r}\left(\sum_{j=k-w}^{\infty}\left|\lambda_{j}\right|\left\|  b_{j}\right\|_{\left(L^p_w,L^q_v \right)_t }\right)^{r} \\
& \lesssim \sum_{k=-\infty}^{\infty}\left(\sum_{j=k-w}^{\infty}\left|\lambda_{j}\right| b^{(k-j) \alpha }\right)^{r} \\
& \lesssim \sum_{k=-\infty}^{\infty}\left(\sum_{j=k-w}^{\infty}\left|\lambda_{j}\right|^{r} b^{(k-j) \alpha r / 2}\right)\left(\sum_{j=k-w}^{\infty} b^{(k-j) \alpha r^{\prime} / 2}\right)^{r/r^{\prime}} \\
& \lesssim \sum_{j=-\infty}^{\infty}\left|\lambda_{j}\right|^{r} \sum_{k=-\infty}^{j+w} b^{(k-j) \alpha r / 2}\\
& \lesssim \sum_{j=-\infty}^{\infty}\left|\lambda_{j}\right|^{r} \\
&\lesssim\|g\|_{\dot{K}{\left(L^p_w,L^q_v \right)_t\left(\mathbb{R}^n \right)}}^{r}.
\end{align*}

Combining these estimates we get 

$$
\|Tg\|_{\dot{K}{\left(L^p_w,L^q_v \right)_t\left(\mathbb{R}^n \right)}} \lesssim\left\|g\right\|_{\dot{K}{\left(L^p_w,L^q_v \right)_t\left(\mathbb{R}^n \right)}}.
$$

Thus, the proof of Theorem  is completed.

\end{proof}
\end{thm}

\section{authors contribution}

Contributions from all authors were equal and significant. The original manuscript was read and approved by all authors.

\section{conflicts of interest}
The authors declare no conflict of interest.

\bibliographystyle{amsplain}

\begin{thebibliography}{99}




  \bibitem{bertrandias}
J. Bertrandias, C. Datry, C. Dupuis,  Unions et intersections d.espaces $L^p$
invariantes par translation
ou convolution, \textit{Ann. Inst. Fourier.} {\bf28}(2) (1978), 53-84.



\bibitem{cowling}
M. Cowling, S. Meda, R. Pasquale,  Riesz potentials and amalgams, \textit{Ann. Inst. Fourier.}
{\bf49}(4) (1999), 1345-1367.

\bibitem{4} D. Cruz-Uribe, A. Fiorenza, J. M. Martell, and C. P\'erez, The boundedness of classical operators on variable $L^p$ spaces, \textit{Ann. Acad. Sci. Fenn. Math.} {\bf 31}(2006), 239–264.
\bibitem{5} D. Cruz-Uribe, A. Fiorenza, and C.J. Neugebauer, The maximal function on variable $L^p$ spaces, \textit{Ann. Acad. Sci. Fenn. Math}. {\bf 28}(2003), 223–238. 


\bibitem{7} L. Diening, Maximal functions on generalized Lebesgue spaces $L^{p(\cdot)}$ , \textit{Math. Inequal.
Appl.} {\bf 7}(2004), 245–253.





\bibitem{fournier}
J. Fournier,  On the Hausdorff-Young theorem for amalgams, \textit{Monatsh. Math.} {\bf95}(2) (1983), 117-135.



\bibitem{herz}
C. Herz, Lipschitz spaces and Bernstein’s theorem on absolutely convergent Fourier transforms, \textit{J. Math. Mech.} {\bf 18} (1968), 283–324.



    \bibitem{archive}
Y. Lu, S. Wang and J. Zhou, Boundedness of some operators on weighted amalgam spaces, \textit{ arXiv}. (2021). https://arxiv.org/abs/2110.01193

\bibitem{23} O. Kov\'acik and J. R\'akosn\'ik, On spaces $L^{p(y)}$ and $W^{k,p(y)}$, \textit{Czechoslovak Math.}
{\bf 41}(116) (1991), 592–618.


\bibitem{32} B. Muckenhoupt, Weighted norm inequalities for the Hardy maximal function,
\textit{Trans. Amer. Math. Soc.} {\bf 165} (1972), 207–226.




	
\bibitem{aims4}
M. Sultan, B. Sultan, A. Khan, T. Abdeljawad, Boundedness of Marcinkiewicz integral operator of variable order in grand Herz-Morrey spaces, \textit{AIMS Math}.  {\bf 8}(9) (2023), 22338-22353. 


\bibitem{ghm4}
B. Sultan, M. Sultan, Boundedness of higher order commutators of Hardy operators on grand Herz-Morrey spaces, \textit{Bull. Sci.math.} { \bf 190}(2024) 103373.


\bibitem{ghm5}
M. Sultan, B. Sultan, A. Hussain, Grand Herz–Morrey Spaces with variable exponent, \textit{Math. Notes.}  {\bf 114} (5) (2023),  957–977.


\bibitem{ghm2}
M. Sultan, B. Sultan, A. Khan, T. Abdeljawad, Boundedness of Marcinkiewicz integral operator of variable order in grand Herz-Morrey spaces, \textit{AIMS Math.} {\bf 8}(9) (2023) 22338-22353.

\bibitem{babargh2}
B. Sultan, M. Sultan, Boundedness of commutators of rough Hardy operators on grand variable Herz spaces, \textit{Forum Math.}  (2023), https://doi .org /10 .1515 /forum -2023 -0152.



\bibitem{babargh4}
B. Sultan, M. Sultan, I. Khan, On Sobolev theorem for higher commutators of fractional integrals in grand variable Herz spaces, \textit{Commun. Nonlinear Sci. Numer. Simul.} {\bf 126} (2023).
\bibitem{ukr}
B. Sultan, M. Sultan,  Sobolev-type Theorem for commutators of Hardy operators in grand Herz spaces, \textit{Ukr. Math. J.} (2024). https://doi.org/10.1007/s11253-024-02381-0



\bibitem{krag}
M. Sultan, B. Sultan, A note on the boundedness of higher order commutators on fractional integrals in grand variable Herz-Morrey spaces, \textit{Kragujev. J. Math.} {\bf 50}(7) (2026), 1063-1080. 

	\bibitem{ineq}
B. Sultan, M. Sultan, A. Khan, T. Abdeljawad, {\ Boundedness of commutators of variable Marcinkiewicz fractional integral operator in grand variable Herz spaces}, \textit{J. Inequal. Appl.} {\bf 2024} (2024), 93.


\bibitem{wiener}
N. Wiener,  On the representation of functions by trigonometrical integrals, \textit{Math. Z.}
{\bf24}(1) (1926), 575-16.

\end{thebibliography}

 \end{document}